\documentstyle{article}
\newtheorem{Definition}{Definition }[section]
\newtheorem{Theorem}[Definition]{Theorem }
\newtheorem{Lemma}[Definition]{Lemma }
\newtheorem{Proposition}[Definition]{Proposition }

\def\proof{\par \noindent{\bf Proof. }\nopagebreak}
\newcommand{\remark}{\vspace{.2cm}\par\noindent{\bf Remark }}
\newcommand{\Remark}[1]{\vspace{.2cm}\par\noindent{\bf Remark }{#1}
\vspace{.2cm}}

\newcommand{\claim}[1]{\vspace{.2cm}\par \noindent{\sc Claim }{#1}
\vspace{.2cm}}
\def\eproof{\nopagebreak\par\nopagebreak\hfill {\Large $\diamond$}\par}

\def\O{{\cal O}}

\def\I{{\cal I}}

\def\O{{\cal O}}
\def\Q{{\bf Q}}

\def\Z{{\bf Z}}
\def\C{{\bf C}}

\def\ra{\rightarrow}

\def\simq{\sim_{\scriptscriptstyle Q}}
\title{Existence of good divisors on Mukai manifolds}
\author{Massimiliano Mella}
\date{}
\begin{document}
\maketitle
\section*{Introduction}
A normal projective variety $X$ is called {\sf Fano} if a multiple
of the \hbox{anticanonical} Weil divisor, $-K_X$, is an ample Cartier
divisor.
The importance of Fano \hbox{varieties} is twofold, from one side they give,
has predicted by Fano \cite{Fa}, \hbox{examples} of non rational varieties
having plurigenera and irregularity all zero (cfr  \cite{Is});
on the other hand they should be the building block of
uniruled variety, indeed recently, Minimal Model Theory predicted that
every
uniruled variety is birational to a fiber space whose general fiber is
a Fano variety with terminal singularities.

The index of a Fano variety $X$ is the number 
$$i(X):=sup\{t\in \Q:
-K_X\equiv tH,\mbox{\rm     for some ample Cartier divisor $H \}$}.$$ It is
known that
$0<i(X)\leq dimX+1$ and if $i(X)\geq dim X$ then $X$ is either an
hyperquadric or a projective space by the Kobayashi--Ochiai criterion,
smooth Fano n-folds
of index $i(X)=n-1$, {\sf del Pezzo n-folds},
have been classified by Fujita \cite{Fu} and terminal Fano n-folds of
index $i(X)>n-2$ have been independently classified by
Campana--Flenner \cite{CF} and Sano \cite{Sa}.

If $X$ has log terminal singularities,
then $Pic(X)$ is torsion free and therefore, the $H$ satisfying
$-K_X\equiv i(X)H$ is uniquely determined and is called the {\sf
fundamental divisor} of $X$.
Mukai announced, \cite{Mu}, the classification of smooth Fano n-folds
$X$
of index $i(X)=n-2$, under the assumption that the linear system $|H|$
contains a smooth divisor.

The main result of this paper is the following

\vspace{.5cm}\noindent {\bf Theorem 1} {\it Let $X$ be a smooth Fano
n-fold of
index $i(X)=n-2$. Then the general element in the fundamental divisor
is smooth.}
\vspace{.5cm}

\noindent Therefore the result of Mukai \cite{Mu} provide a complete
classification of smooth Fano n-folds of index $i(X)=n-2$, {\sf Mukai 
manifolds}.

The ancestors of the theorem, and indeed the lighthouses that
directed
its proof, are Shokurov's proof
for smooth Fano 3-folds, \cite{Sh} and Reid's extension to
canonical Gorenstein 3-folds using the Kawamata's base point free
technique \cite{Re}. This technique was then applied by Wilson in the
case of smooth Fano
\hbox{4-folds} of index 2, \cite{Wi}, afterwards Alexeev, \cite{Al} did it
for
log terminal Fano \hbox{n-folds} of index $i(X)>n-2$ and  recently Prokhorov
used it to prove Theorem 1 in dimension 4 and 5, \cite{Pr1}
\cite{Pr2} \cite{Pr3}. As in Reid's construction we will first
prove the existence of a section with canonical
singularities. To do this we will use Kawamata's
base point free technique and
Kawamata's
notion of centers of log canonical singularities, \cite{Ka1} and his
subadjunction formula for codimension 2 minimal centers \cite{Ka2}.
These tools, together with
Helmke's inductive procedure, \cite{He}, allows to replace difficult
non vanishing arguments by a simple Riemann--Roch calculation. Finally
the Theorem is proved by an inductive argument that lowers the
dimension of $X$.

A natural extension of this problem, motivated by the Minimal Model
Program, should be to ask if for a terminal Fano $X$ of index
$n-2$, with fundamental
divisor $H$ it is true that the general element in $|H|$
has
terminal singularities.
A first, small, step in this direction is the following.

\vspace{.5cm}\noindent {\bf Theorem 2} {\it Let $X$ be a terminal
Gorenstein Fano n-fold of index
n-2. Then the general element in the fundamental divisor $H$ has
canonical singularities.}
\vspace{.5cm}

While working on this subject I had several discussions with M.
Andreatta, who suggested me the direction in which this problem could
be tackled, I would like to express him my deep gratitude, I would
also like to thank A. Corti for valuable comments.

\section{Preliminary results}

We use the standard notation from algebraic geometry.
In particular it is compatible with that of \cite{KMM}
to which we refer constantly, everything is defined over \C.

A \Q-divisor $D$ is an element in $Z_{n-1}(X)\times \Q$,
that is a finite formal sum of prime divisors with rational
coefficients; $D$ is called \Q-Cartier if there is an integer $m$ such
that $mD\in Div(X)$, where $Div(X)$ is the group of Cartier divisors of
$X$. In the following $\equiv$ (respectively $\sim$, $\simq$) will
indicate
numerical  (respectively linear, \Q-linear) equivalence of divisors.
Let $\mu:Y\ra X$ a birational morphism of normal varieties. If
$D$ is a \Q-divisor (\Q-Cartier) then is well defined the strict
transform $\mu_*^{-1}D$ (the pull back $\mu^*D$). For a pair $(X,D)$
of a variety
$X$ and a \Q-divisor $D$, a log resolution is a proper birational
morphism $\mu:Y\ra X$ from a smooth $Y$ such that the union of the
support of $\mu_*^{-1}D$ and of the exceptional locus is a normal
crossing divisor.
\begin{Definition} Let $X$ be a normal variety and $D=\sum_id_iD_i$
an effective
\Q-divisor such that $K_X+D$ is \Q-Cartier. If $\mu:Y\ra X$ is a
log resolution of the pair $(X,D)$, then we can write
$$K_Y+\mu_*^{-1}D=\mu^*(K_X+D)+F$$
with $F=\sum_je_jE_j$ for the exceptional divisors $E_j$. We call
$e_j\in \Q$ the discrepancy coefficient for $E_j$, and regard $-d_i$
as the discrepancy coefficient for $D_i$.

The pair $(X,D)$ is said to have {\sf log canonical} (LC)
(respectively {\sf purely log
terminal} (pLT), {\sf Kawamata log terminal} (KLT)) singularities if
$d_i\leq 1$ (resp. $d_i\leq 1$, $d_i< 1$) and $e_j\geq -1$ (resp.
$e_j>-1$, $e_j>-1$) for any $i,j$ of a log resolution $\mu:Y\ra
X$. In particular if $X$ is smooth at the generic point of $Z$, with
$cod_XZ=a$ and $D$ is a Weil divisor with $mult_ZD=r$, then
$(X,\gamma D)$ is LC for some $\gamma\leq a/r$.
\label{lc}
\end{Definition}
\begin{Definition} A {\sf log-Fano variety} is a pair $(X,\Delta)$
with KLT singularities and such that for some positive integer $m$,
$m(K_X+\Delta)$ is an ample Cartier divisor. The index of a log-Fano
variety $i(X,\Delta):=sup \{t\in \Q: -(K_X+\Delta)\equiv tH$ for some
ample Cartier divisor $H \}$ and the $H$ satisfying
$-(K_X+\Delta)\equiv i(X,\Delta)H$ is called fundamental divisor.
In case $\Delta=0$ we have log terminal Fano variety.
\end{Definition}

We will start recalling
 some results on log
Fano
varieties, essentially due to the Kawamata--Viehweg vanishing theorem.
\begin{Lemma}[\cite{Al}] Let $(X,\Delta)$ be a log-Fano n-fold of index
r, $H$ the fundamental divisor in $X$ and $H^n=d$. Then
\begin{itemize}
\item[-] If $r>n-2$ then $dim |H|=n-2+d(r-n+3)/2>0$
\item[-] If $r=n-2$ and $X$ has canonical Gorenstein singularities,
then $dim |H|=g+n-2\geq n$, where $2g-2=d$, $g\in \Z$, $g\geq 2$.
\end{itemize}
\label{al}
\end{Lemma}

Let us recall the notion and properties of minimal
center of log canonical singularities as introduced in \cite{Ka1}
\begin{Definition}[\cite{Ka1}] Let $X$ be a normal variety and $D=\sum
d_iD_i$ an effective \Q-divisor such that $K_X+D$ is \Q-Cartier.
A subvariety $W$ of $X$ is said to be a {\sf center of log canonical
singularities} for the pair $(X,D)$, if there is a birational morphism
from a normal variety $\mu:Y\ra X$ and a prime divisor $E$ on $Y$ with
the discrepancy coefficient $e\leq -1$ such that $\mu(E)=W$
The set of all the centers of log canonical singularities is denoted
by $CLC(X,D)$, for a point $x\in X$, define $CLC(X,x,D):=\{ W\in
CLC(X,D): x\in W\}$. The union of all the subvarieties in $CLC(X,D)$
is denoted by $LLC(X,D)$.
\end{Definition}
\begin{Theorem}[\cite{Ka1}] Let $X$ be a normal variety and $D$ an
effective \Q-Cartier divisor such that $K_X+D$ is \Q-Cartier. Assume
that $X$ is KLT and $(X,D)$ is LC.
\begin{itemize}
\item[i)]
If $W_1,W_2\in CLC(X,D)$ and $W$ is
an irreducible component of $W_1\cap W_2$, then $W\in CLC(X,D)$. In
particular, if $(X,D)$ is not KLT at a point $x\in X$ then there
exists a unique minimal element of $CLC(X,x,D)$.
\item[ii)]
If $W\in CLC(X,D)$ is a minimal center then $W$ is normal
\item[iii)] Assume that  $D\equiv cL$, with $c<1$, for some ample
Cartier divisor $L$. If $\{x\}\in CLC(X,D)$
is a minimal center then there is a section of $K_X+L$ not vanishing
at $x$.
\end{itemize}
\label{clc}
\end{Theorem}
\begin{Theorem}[\cite{Ka2}] Let $X$ be a normal variety which has only
KLT singularities, $D$ and effective \Q-Cartier divisor such that
$(X,D)$ is LC, and $W$ a minimal element of $CLC(X,D)$. Assume that
$cod W=2$. Then there exist canonically determined effective
\Q-divisors $M_W$ and $D_W$ on $W$ such that $(K_X+D)_{|W}\simq
K_W+M_W+D_W$. If $X$ is affine then there exists an effective
\Q-divisor $M_W^{\prime}$ such that $M_W^{\prime}\simq M_W$ and the
pair $(W,M_W^{\prime}+D_W)$ is KLT.
\label{cod2}
\end{Theorem}
\remark Note that in particular a $cod 2$ minimal center has rational
singularities and if $M_W+D_W\equiv 0$, then $W$ is KLT. In fact It is
enough to choose an open affine covering $\{U_i\}$ of $X$, then for
$V_i=W\cap U_i$ we have $(V_i,M_{V_i}^{\prime}+D_{V_i})$ is KLT and
therefore $V_i$ has rational singularities, \cite{KMM}. Furthermore if
\hbox{$M_W+D_W\equiv 0$} then $M_W\sim D_W\sim \O_W$ and therefore
$M_{V_i}^{\prime}\sim\O_{V_i}$ and these glue together to give that globally $W$ is
KLT.

\begin{Definition}[\cite{He}] Let $X$ be smooth at $x$ and
$(X,D)$ be log canonical at $x$, let $\pi:\tilde{X}\ra X$ the blow up of
$x$.
Following Helmke, the {\sf local discrepancy} of $(X,D)$ is
the rational number
$$
b_x(X,D)=inf\left\{t\left\vert \begin{array}{c} \mbox{\rm
There is a center of log canonical singularity} \\
\mbox{of $(\tilde{X},\pi^*D-(n-1)E+tE)$ contained in $E$}\end{array}
\right.\right\}
$$
\end{Definition}

\claim{ Let $(X,D)$ be LC and $Z\in CLC(X,D)$ with
$x\in Z$ and $Z$, $X$ smooth at $x$, then
$b=b_x(X,D)\leq dim Z$.}

{\it proof of the claim}
Let $\pi:Y\ra X$ the blow up of $x$, with exceptional
divisor $E$ and $Z^{\prime}=\pi^{-1}_*Z\cap E$
since $Z$ is a center of log canonical singularities for $(X,D)$ then
 $\pi^*D$ has multiplicity at least $2cod Z$ along
$Z^{\prime}$. Therefore by definition
$$b\leq -(2cod Z-n+1)+cod Z^{\prime}=dim Z$$.
\eproof

The following inductive procedure due to Helmke (this is a particular
case of his more general Theorem) allows us to decrease the dimension
of a minimal center.

\begin{Proposition}[\cite{He}]
Let $L$ an ample divisor on $X$ and $D$ an effective \Q-divisor with
$D\equiv \gamma L$ for some rational number $0\leq \gamma<1$. Assume
that $X$ is  smooth at
$x$ and $(X,D)$ is log canonical with local discrepancy $b=b_x(X,D)$ at $x$.
Let $Z$ be the minimal center of $CLC(X,x,D)$
 assume that $d=dim Z>0$ and
$Z$, $X$ smooth at $x$.
If
\begin{equation}
L^d\cdot Z>p^d,\hspace{1cm}\mbox{\rm where $p=\frac{b}{1-\gamma}$,}
\label{heq}
\end{equation}
then there is a \Q-divisor
$D_1\equiv \gamma_1L$, with $\gamma<\gamma_1<1$ such that $(X,D_1)$
is log canonical at $x$ with minimal center $Z_1$ properly contained in $Z$ and
$$p_1=\frac{b_1}{1-\gamma_1}<p\hspace{1cm}\mbox{\rm where
$b_1=b_x(X,D_1)$.}$$ \label{He}
\end{Proposition}
We will sketch the proof for reader's convenience.

Step 1. Produce a section $D_0\in |kL_{|Z}|$, for $k\gg0$, with
$mult_x D_0>pk$. This is accomplished by R--R theorem using inequality
(\ref{heq}).

Step 2. Using Serre's vanishing and Bertini Theorem extend $D_0$ to a
section $D^{\prime} \in |kL|$ which
is smooth away from $Z$. Let
$\gamma^{\prime}=sup\{t|(X,D+tD^{\prime})\mbox{\rm
is LC at $x$}\}$, then $D_1=D+\gamma^{\prime}D^{\prime}\equiv
\gamma_1 L$.

Step 3. Use the definition of $p$ and the minimality of $Z$ to prove that
$\gamma_1<1$, and then a straightforward computation gives the assert.

\section{Existence of a canonical section}
For this step we will
use Kawamata's base point free technique, as explained in Reid
\cite{Re}.
Let us start with some lemmas.

\begin{Lemma} Let $X$ be a log terminal Fano n-fold, with $n\geq 3$
and $H$ an ample Cartier divisor with $-K_X\equiv
(n-2)H$ and $G$ a \Q-Cartier divisor with $(X,G)$ LC. Assume that
$Z\in CLC(X,G)$ is a minimal center and $G\equiv
\gamma H$, with $\gamma< cod Z-1$
then there is a section
of $H$ not vanishing identically on $Z$.
\label{KV}
\end{Lemma}

\proof
We proceed as in \cite[Prop 2.3]{Ka1}. 
Let $M\in |mH|$, for $m\gg 0$, be a general member among 
Cartier divisors containing $Z$, let \hbox{$G_1=
(1-\epsilon_1) (G+\epsilon_2 M)$,} for $\epsilon_i\ll 1/m$, 
 then $G_1\equiv \gamma_1 H$, with 
$\gamma_1<cod Z-1$. Furthermore we may assume that 
$(X,G_1)$ is LC  and
$Z$ is an isolated element of $LLC(X,G_1)$.
Let $\mu:Y\ra X$ a log resolution of $(X,G_1)$, then
$$K_Y+E-A+F=\mu^*(K_X+G_1),$$
where $\mu(E)=Z$, $A$ is an integral \hbox{$\mu$-exceptional}
divisor and $\lfloor F\rfloor=0$.
Let 
$$ N(t):= -K_Y-E-F+A+\mu^*(tH),$$
then $N(t)\equiv \mu^*(t+(n-2)-\gamma_1)H$ and $N(t)$ is nef and big 
whenever
$t+(n-2)-\gamma_1> 0$, hence by hypothesis this is true whenever $t\geq
-n+1+cod Z$.
Thus K--V vanishing yields
\begin{equation}
H^i(Y,\mu^*(tH)-E+A)=0\hspace{.7cm}
H^i(E,(\mu^*(tH)+A)_{|E_0})=0
\label{van}
\end{equation}
for $i>0$ and $t\geq -n+1+cod Z$,
and consequently
$$
 H^0(Y,\mu^*H+A)\ra H^0(E,\mu^*H+A)\ra 0;
$$
since $A$ is effective and $\mu$-exceptional, then any section
 in $H^0(Y,\mu^*H+A)$, not vanishing on $E$, pushes forward to
give a section of $H$ not vanishing on $Z$.
To conclude the proof it is, therefore, enough to prove that
$h^0(E,N(1))>0$. Let
$p(t)=\chi(E,N(t))$,
then by equation (\ref{van}), 
$p(0)\geq 0$ and $p(t)=0$ for $0>t\geq -n+1+cod Z=-dim
Z+1$. Since $deg p(t)= dim Z$ and $p(t)>0$ for $t\gg 0$
then $h^0(E,N(1))=p(1)>0$.
\eproof

The above lemma allows us, essentially, to treat minimal centers of
codimension $\geq 3$. In the next couple of lemmas we will treat
codimension 2 minimal centers.

\begin{Lemma} Let $X$ be a log terminal Fano n-fold, with $n\geq 3$,
and $H$ an ample Cartier divisor with $-K_X\equiv
(n-2)H$, let $L\sim (n-1)H$ and $D$ a \Q-Cartier divisor with $(X,D)$
LC. Assume that $D\equiv H$ and $Z\in CLC(X,D)$ a cod 2 minimal
center, then for $k\gg 0$ and $\delta\geq 0$
$$h^0(Z,kL_{|Z})\geq
\frac{(n-1)^{n-2}}{(n-2)!}
k^{n-2}+\frac{(n-3+\delta)(n-1)^{n-3}}{2(n-3)!}k^{n-3}+\mbox{\rm
lower terms in k}.$$
Furthermore, keeping the notation of Theorem \ref{cod2}, if
$M_Z+D_Z\not\equiv 0$ then $\delta>0$.
\label{K2}
\end{Lemma}
\proof
By Theorem \ref{cod2} there are effective \Q-divisors $M_Z$
and $D_Z$ such that
$$ -(n-3)H_{|Z}\equiv (K_X+D)_{|Z}\simq K_Z+M_Z+D_Z.$$

 Let $f:Y\ra Z$ a
log resolution of $(Z,M_Z+D_Z)$
then $$K_Y+\Delta=f^*(K_Z+M_Z+D_Z)+\sum e_iE_i,$$
 where $\Delta=f^{-1}_*(M_Z+D_Z)$ is
effective and the $E_i$ are $f$-exceptional.
In particular
$$-K_Y\cdot (f^*H_{|Z})^{n-3}=-(K_X+D)_{|Z}\cdot H^{n-3}_{|Z}+\Delta\cdot
f^*H_{|Z}^{n-3}\geq n-3+\delta$$
Since $Z$ has
rational singularities and $L$ is ample then,  for $k\gg 0$,
$$h^0(Z,kL_{|Z})=\chi(Z,kL_{|Z})=\chi(Y,kf^*L_{|Z}).$$ 
Hence by R--R
formula $$h^0(Z,kL_{|Z})\geq
\frac{(n-1)^{n-2}}{(n-2)!}
k^{n-2}+\frac{(n-3+\delta)(n-1)^{n-3}}{2(n-3)!}k^{n-3}+\mbox{\rm
lower terms in k}.$$
\eproof

\begin{Lemma} Let $X$ be a log terminal Gorenstein Fano n-fold, with
$n\geq 3$, and $H$ an ample Cartier divisor with $-K_X\sim
(n-2)H$, let $L\sim (n-1)H$ and $D$ a \Q-Cartier divisor with $(X,D)$
LC. Assume that $D\equiv H$ and $Z\in CLC(X,D)$ a cod 2 minimal
center with
$Z\not\subset Sing(X)$. Then there exists
a section of $H\sim K+L$ not vanishing identically on $Z$.
\label{KH}
\end{Lemma}
\proof
By Theorem \ref{cod2} there are effective \Q-divisors $M_Z$
and $D_Z$ such that
$$ -(n-3)H_{|Z}\equiv (K_X+D)_{|Z}\simq K_Z+M_Z+D_Z.$$

If $n=3$ then $Z$ is a smooth curve, by Theorem \ref{clc}, and 
$g(Z)\leq 0$, thus
$h^0(Z,H)>0$; if $n>3$ and $(Z,M_Z+D_Z)$ is KLT
then $(Z,M_Z+D_Z)$
is a log-Fano variety of index $i(Z,M_Z+D_Z)=dimZ-1$, therefore by Lemma
\ref{al}, $h^0(Z,H)>0$.
As in the proof of Lemma \ref{KV} let us replace $D$ with $D_1$ such that
$(X,D_1)$ is LC, $Z$ is isolated in
$LLC(X,D_1)$ and $D_1\equiv \gamma_1H$, for 
$\gamma_1<1+\epsilon$, with $\epsilon\ll 1$.
Let $\mu:Y\ra X$ a log resolution of $(X,D_1)$ with
$K_Y+E-A+F=\mu^*(K_X+D_1)$, where $f(E)=Z$, $A$ is an integral 
\hbox{$\mu$-exceptional}
divisor and $\lfloor F\rfloor=0$. 
Let
$N(t):=-K_Y-E-F+A+\mu^*(tH)$, then
$N(1)\equiv \mu^*(1+(n-2)-\gamma_1)H$ is nef and big and
consequently
$$
 H^0(Y,\mu^*H+ A)\ra H^0(E,\mu^*H+ A)\ra 0.
$$
Therefore the sections in $H^0(Z,H)$ extends to sections of $H^0(X,H)$ not
vanishing identically on $Z$.

\noindent By the remark after Theorem
\ref{cod2} we can, therefore assume
that $M_Z+D_Z\not\equiv 0$.
Fix a smooth point $x\in Z$ outside of $Sing(X)$, such that $Z$ is 
the minimal element of $CLC(X,x,D)$.
Let us mimic Helmke's arguments;
in the notation of Proposition \ref{He},
$\gamma=1/(n-1)$ and
$$p=\frac{b}{1-\gamma}\leq \frac{(n-1)(n-2)}{n-2}\leq n-1.$$
The first
step is accomplished using Lemma \ref{K2}; in fact
$$h^0(Z,\O_Z/\I^k_{Z,x})=\frac{1}{(n-2)!}k^{n-2}+
\frac{(n-3)}{2(n-3)!}k^{n-3}+\mbox{\rm
lower terms in k},$$
therefore by Lemma \ref{K2} there exists a section $D^{\prime}\in|kL_{|Z}|$,
for $k\gg 0$, such that \hbox{$mult_xD^{\prime}>pk$}.
 It is now enough to follow word by word Helmke's
arguments to conclude that there is a
\Q-divisor $D_1\equiv \gamma_1L$, with $\gamma<\gamma_1<1$ such that
$(X,D_1)$ is log
canonical at $x$, with minimal center $Z_1\ni x$
properly contained in $Z$.
Since $x\in Z_1$ and $p_1<p\leq n-1$
then $Z_1\not\subset Sing(X)$ and we can choose a smooth point
$x_1\in Z_1$ and apply directly Proposition \ref{He} to $(X,D_1)$ and $x_1$.
Inductively the dimension of the minimal center is lowered and we find
a divisor $D_l\equiv \gamma_lL$, with $c_l<1$,
which has zero dimensional minimal center. Conclude by
 Theorem \ref{clc} iii) that there exists a section of $H\sim
K_X+L$ not vanishing on $Z$.
\eproof

We will need the forthcoming lemma only in the next section,
to be able
to apply an inductive procedure on the Fano variety, but we place it
here since the flavor and the proof are close to the previous one.

\begin{Lemma}  Let $X$ be a log terminal Gorenstein Fano n-fold, with
$n>3$,	and $H$ an ample Cartier divisor with $-K_X\sim
(n-2)H$, let $L\sim (n-1)H$ and $D$ a \Q-Cartier divisor with $(X,D)$
LC. Assume that
$D\equiv  2H$, $Z\in
CLC(X,D)$ is a codimension 3 minimal center not contained in
$Sing(X)$. Furthermore assume that
there exist
$S\in|H|$ and an effective \Q-divisor $D_S$
such that $(S,D_S)$ satisfy the
hypothesis of Theorem \ref{cod2} and Lemma \ref{K2}.  Then there is
a section of $H\sim K_X+L$ not vanishing identically on $Z$.
\label{KH2}
\end{Lemma}
\proof Let us simply sketch the proof since it is similar to that of
Lemma \ref{KH}.
By Theorem \ref{cod2} there exist effective \Q-divisors $M_Z$ and
$D_Z$ such that
$$-(n-4)H_{|Z}\equiv K_Z+M_Z+D_Z. $$
Let us, again, replace $D$ with $D_1$ such that
$(X,D_1)$ is LC, $Z$ is isolated in
$LLC(X,D_1)$ and $D_1\equiv \gamma_1H$, for 
$\gamma_1<2+\epsilon$, with $\epsilon\ll 1$.
Let $\mu:Y\ra X$ a log resolution of $(X,D_1)$ with
$K_Y+E-A+F=\mu^*(K_X+D_1)$, where $f(E)=Z$, $A$ is an integral 
$\mu$-exceptional
divisor and $\lfloor F\rfloor=0$. 

Let
\hbox{$N(t):=-K_Y-E-F+A+\mu^*(tH)$,} then
$N(1)\equiv \mu^*(1+(n-2)-\gamma_1)H$ is nef and big

\noindent If $n=4$ then $Z$ is a smooth curve of non positive 
genus, therefore
$h^0(Z,H)>0$; if $n\geq 5$ and
$(Z,M_Z+D_Z)$ is KLT then it is a log-Fano variety of index $i(Z,M_Z+D_Z)=
dim Z-1$, therefore as above the sections in $H^0(Z,H)$ extends to
sections of $H^0(X,H)$ not vanishing identically on $Z$.

Again we can assume that $M_Z+D_Z\not\equiv 0$, and choose
a smooth point in $Z$ with $x\not\in Sing(X)$, in Helmke's notations,
$\gamma=2/(n-1)<1$ and
$$p\leq \frac{(n-1)(n-3)}{n-1-2}\leq n-1$$
and by Lemma \ref{K2} for $k\gg 0$
there is a section $D^{\prime}\in |kL_{|Z}|$ with $mult_x D^{\prime}>pk$;
then
conclude as in Lemma \ref{KH}
\eproof

\begin{Proposition} Let $X$ be a log terminal Gorenstein Fano n-fold
and $H$ an ample Cartier divisor with $-K_X\sim
(n-2)H$. Assume that
$codSing(X)>2$ and $n\geq 3$.
Then the general element in $|H|$ has at worst canonical singularities.
\label{alto}
\end{Proposition}

\proof By Lemma \ref{al} we know that $dim |H|\geq 1$. Let $S\in |H|$
a generic element and assume that $S$
has worse than canonical singularities.
Since both $H$ and $K_X$ are Cartier divisors then $(X,S)$ is not pLT,
that is there exists $\gamma\leq 1$ such that $(X,\gamma S)$,
is LC with $Z$ a minimal
center in $CLC(X,\gamma S)$, and by Bertini theorem
$Z\subset Bsl|H|$.
We will derive	a contradiction ,
producing a section of $|H|$ not vanishing identically on $Z$.

\noindent If either $dim Z\leq n-3$ or $dim Z=n-2$ and $\gamma<1$ then
apply Lemma \ref{KV}.

\noindent If $dim Z=n-2$ and $\gamma=1$,
by hypothesis $Z\not\subset Sing(X)$ hence apply
Lemma \ref{KH}.

To conclude we have to exclude the case $dim Z=n-1$.
Assume that
$|H|$ has a fixed component $F$, by \cite[Prop 3.2]{Al}
$F$ must have multiplicity
1, that is $\gamma=1$.	Since $H$ is connected and
movable then $S$ must be singular along a codimension 2 set $Z\subset
F$, therefore $F$ is not minimal in $CLC(X,S)$, see Definition \ref{lc}.
\eproof

\Remark{
In particular the above argument shows that $H$ is smooth in codimension
1 and there are not fixed component.}

\proof (of Theorem 2) By Lemma \ref{al}, $h^0(X,H)\geq 2$; furthermore
terminal singularities are
smooth in codimension 2. It is, therefore  enough to apply Proposition
\ref{alto}.
\eproof

\remark It is not true, in general, that terminal Fano $X$ of index
$i(X)=n-2$ are Gorenstein; consider a terminal Fano
3-fold
with an Enriques surface as section of the fundamental divisor, this
varieties are studied by Conte--Murre \cite{CM}. In this case $X$ has
8 singular points, which are cones over the Veronese surface, and
$X$ is 2-Gorenstein; nevertheless $H$ has a smooth (terminal) section.

\section{Proof of the main Theorem}

By a direct calculation, for instance Lemma \ref{al},  $h^0(X,H)\geq
n$ therefore  by
Proposition
\ref{alto} there exists a section $S\in |H|$ with canonical
singularities. Our aim is to apply inductively Proposition
\ref{alto}, to do this we have to prove that $S$ is smooth in
codimension 2.
Assume the contrary, in particular $S$ is not terminal and there is a
center $Z\subset Bsl|H|$ of canonical singularities in $S$ with $dim
Z=n-3$.

{\sc Case 1} Assume that all
sections of $|H|$ are singular at $Z$, let
$H_i\in |H|$ generic elements and \hbox{$D=1/2(H_1+H_2)$}.

\claim{$(X,\gamma D)$ is log canonical for some
$\gamma\leq 3/2$ with a minimal
center $W\subseteq Z$ of codimension $\geq 3$.}

Observe that by the claim we can apply Lemma \ref{KV}, to produce a
section of $|H|$ not vanishing on $Z$ and derive in this way a
contradiction.

\proof(of the claim) Let $f:Y\ra X$ the blow up of $Z$
let $f^*S=S^{\prime}+ rE$, since $X$ is
smooth at the generic point of $Z$ then
$K_Y=f^*K_X+2E$. By adjunction formula
$$K_{S^{\prime}}=(K_Y+S^{\prime})_{|S^{\prime}}=f^*K_{S}+(2-r)
E_{|S^{\prime}},$$
since $S$ is canonical and is singular at
$Z$ then $r=2$.
$|H|$ has not fixed components and its general
element
is smooth in codimension 1 therefore for some $\gamma\leq
3/2$, $(X,\gamma D)$ is log canonical with a minimal center
$W\subseteq Z$ of codimension $\geq 3$.
\eproof

{\sc Case 2} Assume that there are infinitely many such codimension 3
components $Z_i\subset Bsl|H|$ centers of canonical singularities for
$H_i\in |H|$. Let $H_1$ a generic element in $|H|$,
we can
assume that $H_1$ is singular along $Z_1$, with $Z_1\subset Bsl|H|$
and
$Z_1\not\subset Sing(X)$.
Let $D=1/2(H_1+H_2)$, with $H_2\in |H|$, a general element;
by construction $(X,\gamma D)$ is log canonical for some $\gamma \leq
2$ with a
minimal center $Z$ of codimension$\geq 3$. If either $\gamma <2$ or
$codZ>3$ then conclude by Lemma \ref{KV}.
Assume that $\gamma=2$ and $cod Z=3$, we can assume without loss of
generality that $Z=Z_1$, let $S\in |H|$ a generic element smooth
at the generic point of $Z$ and $D_S=H_{1|S}$, then
$(S,D_S)$ and
$Z$ satisfy the hypothesis of Theorem \ref{cod2} and Lemma \ref{K2},
thus, we derive a contradiction by
Lemma \ref{KH2} if $n\geq
4$.

At each inductive step we loose only one section of $|H|$, therefore
$|H_{|S}|$ is always movable; furthermore
by K--V vanishing theorem
$$H^0(X,H)\ra H^0(S,H_{|S})\ra 0,$$
hence it is possible to
study the singularities of $S$ trough the linear system $H_{|S}$.
To carry on induction in Case 1 we need the following
\claim{ Let $S_0=X$ and $S_j=X\cap H_1\cap\ldots\cap H_j$, for $H_i\in
|H|$
general elements. Assume that $S_j$ has canonical singularities and is
singular at
$Z_j$, with $cod_XZ_j=j+2$. If $X$ is smooth at $Z_j$ then $S_{j-1}$
is smooth at $Z_j$.}
\proof(of the claim) We will prove it by induction on $j$. If $j=1$
then it follows by hypothesis. Let $f:Y\ra X$ the blow
up of $Z_j$, with $f^*H_i=H_i^{\prime}+ rE$, since $X$ is
smooth at the generic point of $Z_j$ then
$K_Y=f^*K_X+(j+1)E$. By adjunction formula
$$K_{S_j^{\prime}}=(K_Y+\sum_i
H_i^{\prime})_{|S_j^{\prime}}=f^*K_{S_j}+(j+1-jr) E_{|S_j^{\prime}},$$
where $S_j^{\prime}=Y\cap H_1^{\prime}\cap\ldots\cap H_j^{\prime}$.
Since $S_j$ has canonical singularities then \hbox{$j+1-jr\geq 0$,} and
consequently $r=1$, that is the generic element in $|H|$ is
smooth at $Z_j$. On the other hand $Z_j=S_1\cap H_2\cap\ldots\cap
H_j$, where $S_1\in |H|$ is a general element smooth at $Z_j$, and
$cod_{S_1}Z_j=j-1$, therefore by induction hypothesis $S_{j-1}$ is
smooth at $Z_j$.
\eproof

 By the inductive process we are reduced to
a canonical Gorenstein 3-fold smooth in codimension 2,
$S_3=X\cap(\bigcap_{i=1}^{n-3}H_i)$ with a line bundle
$H_3=H_{|S_3}$ satisfying the following conditions:
\begin{itemize}
\item[-] $h^0(S_3,H_3)\geq 3$
\nopagebreak\item[-] $dim Bsl|H_3|= 1$
\end{itemize}
Let $H_1\in |H_3|$ a general element and $B$ a curve contained in the
base locus of $|H_3|$. Assume, without loss of generalities that
$x_1\in B\cap H_1$ is such that
$x_1\not\in Sing(S_3)$ and $H_1$ singular at $x_1$.
Let $A=\{ M\in|H|$ $|$ $M$  is singular
at $x_1 \}$ since $h^0(S_3,H_3)\geq 3$ and $dim Bsl|H_3|=1$
then $dim A\geq 1$. Let $H_i\in A$, for
$i=1,2$ be general elements
and $D=1/2(D_1+D_2)$. Then $(X,\gamma D)$ is log canonical for some
$\gamma\leq 3/2$ with zero dimensional	minimal center
thus Lemma \ref{KV} apply to derive a contradiction and the theorem is proved.
\eproof
\begin{flushleft}
Massimiliano Mella\\
Universit{\`a} degli Studi di Trento\\
Dipartimento di Matematica\\
I38050 Povo (TN), Italia\\
e-mail: mella@science.unitn.it
\end{flushleft}
\end{document}